\documentclass[12pt]{article}
\usepackage{amsmath,amssymb,amsthm,amscd,a4wide}

\begin{document}

\newcommand{\End}{{\rm{End}\ts}}
\newcommand{\Hom}{{\rm{Hom}}}
\newcommand{\ch}{{\rm{ch}\ts}}
\newcommand{\non}{\nonumber}
\newcommand{\wt}{\widetilde}
\newcommand{\wh}{\widehat}
\newcommand{\ot}{\otimes}
\newcommand{\la}{\lambda}
\newcommand{\La}{\Lambda}
\newcommand{\De}{\Delta}
\newcommand{\al}{\alpha}
\newcommand{\be}{\beta}
\newcommand{\ga}{\gamma}
\newcommand{\ve}{\varepsilon^{}}
\newcommand{\ep}{\epsilon}
\newcommand{\ka}{\kappa}
\newcommand{\vk}{\varkappa}
\newcommand{\si}{\sigma}
\newcommand{\vp}{\varphi}
\newcommand{\de}{\delta^{}}
\newcommand{\ze}{\zeta}
\newcommand{\om}{\omega}
\newcommand{\hra}{\hookrightarrow}
\newcommand{\ee}{\ep^{}}
\newcommand{\su}{s^{}}
\newcommand{\ts}{\,}
\newcommand{\vac}{\mathbf{1}}
\newcommand{\di}{\partial}
\newcommand{\qin}{q^{-1}}
\newcommand{\tss}{\hspace{1pt}}
\newcommand{\Sr}{ {\rm S}}
\newcommand{\U}{ {\rm U}}
\newcommand{\X}{ {\rm X}}
\newcommand{\Y}{ {\rm Y}}
\newcommand{\Z}{{\rm Z}}
\newcommand{\AAb}{\mathbb{A}\tss}
\newcommand{\CC}{\mathbb{C}\tss}
\newcommand{\QQ}{\mathbb{Q}\tss}
\newcommand{\SSb}{\mathbb{S}\tss}
\newcommand{\ZZ}{\mathbb{Z}\tss}
\newcommand{\Ac}{\mathcal{A}}
\newcommand{\Lc}{\mathcal{L}}
\newcommand{\Mc}{\mathcal{M}}
\newcommand{\Pc}{\mathcal{P}}
\newcommand{\Qc}{\mathcal{Q}}
\newcommand{\Tc}{\mathcal{T}}
\newcommand{\Sc}{\mathcal{S}}
\newcommand{\Bc}{\mathcal{B}}
\newcommand{\Dc}{\mathcal{D}}
\newcommand{\Ec}{\mathcal{E}}
\newcommand{\Fc}{\mathcal{F}}
\newcommand{\Hc}{\mathcal{H}}
\newcommand{\Uc}{\mathcal{U}}
\newcommand{\Wc}{\mathcal{W}}
\newcommand{\Zc}{\mathcal{Z}}
\newcommand{\Ar}{{\rm A}}
\newcommand{\Ir}{{\rm I}}
\newcommand{\Gr}{{\rm G}}
\newcommand{\Or}{{\rm O}}
\newcommand{\Spr}{{\rm Sp}}
\newcommand{\Wr}{{\rm W}}
\newcommand{\Zr}{{\rm Z}}
\newcommand{\gl}{\mathfrak{gl}}
\newcommand{\Pf}{{\rm Pf}\ts}
\newcommand{\oa}{\mathfrak{o}}
\newcommand{\spa}{\mathfrak{sp}}
\newcommand{\g}{\mathfrak{g}}
\newcommand{\h}{\mathfrak h}
\newcommand{\n}{\mathfrak n}
\newcommand{\z}{\mathfrak{z}}
\newcommand{\Zgot}{\mathfrak{Z}}
\newcommand{\p}{\mathfrak{p}}
\newcommand{\sll}{\mathfrak{sl}}
\newcommand{\agot}{\mathfrak{a}}
\newcommand{\qdet}{ {\rm qdet}\ts}
\newcommand{\Ber}{ {\rm Ber}\ts}
\newcommand{\Mat}{{\rm{Mat}}}
\newcommand{\HC}{ {\mathcal HC}}
\newcommand{\cdet}{ {\rm cdet}}
\newcommand{\tr}{ {\rm tr}}
\newcommand{\gr}{{\rm gr}\ts}
\newcommand{\str}{ {\rm str}}
\newcommand{\loc}{{\rm loc}}
\newcommand{\sgn}{ {\rm sgn}\ts}
\newcommand{\hb}{\mathbf{h}}
\newcommand{\ba}{\bar{a}}
\newcommand{\bb}{\bar{b}}
\newcommand{\bi}{\bar{\imath}}
\newcommand{\bj}{\bar{\jmath}}
\newcommand{\bk}{\bar{k}}
\newcommand{\bl}{\bar{l}}
\newcommand{\Sym}{\mathfrak S}
\newcommand{\fand}{\quad\text{and}\quad}
\newcommand{\Fand}{\qquad\text{and}\qquad}
\newcommand{\For}{\qquad\text{or}\qquad}
\newcommand{\vt}{{\tss|\hspace{-1.5pt}|\tss}}

\renewcommand{\theequation}{\arabic{section}.\arabic{equation}}

\newtheorem{thm}{Theorem}[section]
\newtheorem{lem}[thm]{Lemma}
\newtheorem{prop}[thm]{Proposition}
\newtheorem{cor}[thm]{Corollary}
\newtheorem{conj}[thm]{Conjecture}
\newtheorem*{mthm}{Main Theorem}

\theoremstyle{definition}
\newtheorem{defin}[thm]{Definition}

\theoremstyle{remark}
\newtheorem{remark}[thm]{Remark}
\newtheorem{example}[thm]{Example}

\newcommand{\bth}{\begin{thm}}
\renewcommand{\eth}{\end{thm}}
\newcommand{\bpr}{\begin{prop}}
\newcommand{\epr}{\end{prop}}
\newcommand{\ble}{\begin{lem}}
\newcommand{\ele}{\end{lem}}
\newcommand{\bco}{\begin{cor}}
\newcommand{\eco}{\end{cor}}
\newcommand{\bde}{\begin{defin}}
\newcommand{\ede}{\end{defin}}
\newcommand{\bex}{\begin{example}}
\newcommand{\eex}{\end{example}}
\newcommand{\bre}{\begin{remark}}
\newcommand{\ere}{\end{remark}}
\newcommand{\bcj}{\begin{conj}}
\newcommand{\ecj}{\end{conj}}

\newcommand{\bal}{\begin{aligned}}
\newcommand{\eal}{\end{aligned}}
\newcommand{\beq}{\begin{equation}}
\newcommand{\eeq}{\end{equation}}
\newcommand{\ben}{\begin{equation*}}
\newcommand{\een}{\end{equation*}}

\newcommand{\bpf}{\begin{proof}}
\newcommand{\epf}{\end{proof}}

\def\beql#1{\begin{equation}\label{#1}}

\title{\Large\bf Pfaffian-type Sugawara operators}

\author{A. I. Molev\\[2em]
School of Mathematics and Statistics\\
University of Sydney,
NSW 2006, Australia\\
alexander.molev@sydney.edu.au}

\date{} 
\maketitle


\begin{abstract}
We show that the Pfaffian of a generator matrix for the affine
Kac--Moody algebra $\wh\oa_{2n}$ is a Segal--Sugawara vector.
Together with our earlier
construction involving the symmetrizer in the Brauer algebra, this gives
a complete set of Segal--Sugawara vectors in type $D$.
\end{abstract}


%

\section{Introduction}
\label{sec:int}
\setcounter{equation}{0}

For each affine Kac--Moody algebra $\wh\g$ associated with
a simple Lie algebra $\g$, the corresponding vacuum module $V(\g)$ at the critical level
is a vertex algebra.
The structure of the center $\z(\wh\g)$ of this vertex algebra
was described by
a remarkable theorem of Feigin
and Frenkel in \cite{ff:ak}, which states that $\z(\wh\g)$
is the algebra of polynomials in infinitely many variables
which are associated with generators of the algebra of $\g$-invariants
in the symmetric algebra $\Sr(\g)$. For a detailed proof of the theorem,
its extensions and significance for the representation theory of
the affine Kac--Moody algebras see \cite{f:lc}.

In a recent paper \cite{m:ff}
we used the symmetrizer in the Brauer algebra to construct
families of elements of $\z(\wh\g)$
(Segal--Sugawara vectors) for the Lie algebras $\g$ of types $B$, $C$ and $D$
in an explicit form.
In types $B$ and $C$ they include
complete sets of Segal--Sugawara vectors generating the center
$\z(\wh\g)$, while
in type $D$ one vector in the complete set of \cite{m:ff}
was not given explicitly. The aim of this note is to produce
this Segal--Sugawara vector in $\z(\wh\oa_{2n})$ which is associated with the Pfaffian
invariant in $\Sr(\oa_{2n})$.

Simple explicit formulas for generators of the Feigin--Frenkel
center $\z(\wh\gl_n)$ were given recently in
\cite{cm:ho}, \cite{ct:qs}, following Talalaev's construction
of higher Gaudin Hamiltonians \cite{t:qg}. So, together with 
the results of \cite{m:ff} 
we get a construction of generators of the Feigin--Frenkel
centers for all classical types.

\section{Pfaffian-type generators}
\label{sec:pg}
\setcounter{equation}{0}

Denote by $E_{ij}$, $1\leqslant i,j\leqslant 2n$,
the standard basis vectors of the Lie
algebra $\gl_{2n}$.
Introduce the elements $F_{ij}$
of $\gl_{2n}$ by the formulas
\beql{fijelem}
F_{ij}=E_{ij}-E_{ji}.
\eeq
The Lie subalgebra
of $\gl_{2n}$ spanned by the elements $F_{ij}$ is isomorphic to
the even orthogonal Lie algebra $\oa_{2n}$.
The elements of $\oa_{2n}$ are skew-symmetric matrices. Introduce
the standard normalized invariant bilinear form on $\oa_{2n}$ by
\ben
\langle X,Y\rangle=\frac{1}{2}\ts\tr\tss XY,\qquad X,Y\in\oa_{2n}.
\een

Now consider the affine
Kac--Moody algebra $\wh\oa_{2n}=\oa_{2n}\tss[t,t^{-1}]\oplus\CC K$
and set $X[r]=Xt^r$ for any $r\in\ZZ$ and $X\in\oa_{2n}$. The element $K$
is central in $\wh\oa_{2n}$ and
\ben
\big[X[r],Y[s]\big]=[X,Y][r+s]+r\ts\de_{r,-s}\langle X,Y\rangle\ts K.
\een
Therefore, for the generators we have
\ben
\bal
\big[F_{ij}[r],F_{kl}[s]\big]&=\de_{kj}\ts F_{il}[r+s]-\de_{il}\ts F_{kj}[r+s]
-\de_{ki}\ts F_{jl}[r+s]+\de_{jl}\ts F_{ki}[r+s]\\
{}&+r\tss\de_{r,-s}\big(\de_{kj}\ts\de_{il}-\de_{ki}\ts\de_{jl}\big)\ts K.
\eal
\een
The {\it vacuum module at the critical level} $V(\wh\oa_{2n})$ can be defined
as the quotient of the
universal enveloping algebra $\U(\wh\oa_{2n})$ by the left
ideal generated by $\oa_{2n}[t]$ and $K+2n-2$ (note that the dual Coxeter number in type $D_n$
is $h^{\vee}=2n-2$).
The Feigin--Frenkel center $\z(\wh\oa_{2n})$ is defined by
\ben
\z(\wh\oa_{2n})=\{v\in V(\wh\oa_{2n})\ |\ \oa_{2n}[t]\ts v=0\}.
\een
Any element of $\z(\wh\oa_{2n})$ is called a {\it Segal--Sugawara vector}.
A {\it complete set of Segal--Sugawara vectors\/}
$\phi_{2\ts 2},\phi_{4\ts 4},\dots,\phi_{2n-2\ts 2n-2},\phi'_n$
was produced in \cite{m:ff}, where all of them,
except for $\phi'_n$, were given explicitly. We will produce $\phi'_n$
in Theorem~\ref{thm:pfaff} below.

Combine the generators $F_{ij}[-1]$ into the skew-symmetric
matrix $F[-1]=\big[F_{ij}[-1]\big]$
and define its {\it Pfaffian} by
\ben
\Pf F[-1]=\frac{1}{2^nn!}\sum_{\si\in\Sym_{2n}}\sgn\si\cdot
F_{\si(1)\ts\si(2)}[-1]\dots F_{\si(2n-1)\ts\si(2n)}[-1].
\een
Note that the elements $F_{ij}[-1]$ and $F_{kl}[-1]$ of $\wh\oa_{2n}$ commute,
if the indices $i,j,k,l$ are distinct. Therefore, we can write the formula
for the Pfaffian in the form
\beql{pf}
\Pf F[-1]=\sum_{\si}\sgn\si\cdot
F_{\si(1)\ts\si(2)}[-1]\dots F_{\si(2n-1)\ts\si(2n)}[-1],
\eeq
summed over the elements $\si$ of the subset $\Bc_{2n}\subset \Sym_{2n}$
which consists of the permutations with the properties
$\si(2k-1)<\si(2k)$ for all $k=1,\dots,n$ and $\si(1)<\si(3)<\dots<\si(2n-1)$.

\bth\label{thm:pfaff}
The element $\phi'_n=\Pf F[-1]$ is a Segal--Sugawara vector for $\wh\oa_{2n}$.
\eth

\bpf
We need to show that $\oa_{2n}[t]\ts \phi'_n=0$ in the vacuum module $V(\wh\oa_{2n})$.
It suffices to verify that for all $i,j$,
\beql{annih}
F_{ij}[0]\ts \Pf F[-1]=F_{ij}[1]\ts \Pf F[-1]=0.
\eeq
Note that for any permutation $\pi\in\Sym_{2n}$ the mapping
\ben
F_{ij}[r]\mapsto F_{\pi(i)\ts\pi(j)}[r],\qquad K\mapsto K
\een
defines an automorphism of the Lie algebra $\wh\oa_{2n}$. Moreover,
the image of $\Pf F[-1]$ under its extension to $\U(\wh\oa_{2n})$ coincides
with $\sgn\pi\cdot\Pf F[-1]$. Hence, it is enough to verify \eqref{annih}
for $i=1$ and $j=2$.

Observe that $F_{12}[0]$ commutes with all summands in \eqref{pf} with
$\si(1)=1$ and $\si(2)=2$. Suppose now that $\si\in\Bc_{2n}$
is such that $\si(2)>2$. Then $\si(3)=2$ and $\si(4)>2$.
In $V(\wh\oa_{2n})$ we have
\ben
\bal
F_{12}[0]\ts F_{1\ts \si(2)}[-1]\tss F_{2\ts \si(4)}[-1]&\dots F_{\si(2n-1)\ts\si(2n)}[-1]\\
{}={}&{}-F_{2\ts \si(2)}[-1]\tss F_{2\ts \si(4)}[-1]\dots F_{\si(2n-1)\ts\si(2n)}[-1]\\
{}&{}+F_{1\ts \si(2)}[-1]\tss F_{1\ts \si(4)}[-1]\dots F_{\si(2n-1)\ts\si(2n)}[-1].
\eal
\een
Set $i=\si(2)$ and $j=\si(4)$. Note that the permutation $\si'=\si\ts (2\ts 4)$ also belongs
to the subset $\Bc_{2n}$, and $\sgn\si'=-\sgn\si$. We have
\begin{multline}
-F_{2\ts i}[-1]\tss F_{2\ts j}[-1]+F_{1\ts i}[-1]\tss F_{1\ts j}[-1]
+F_{2\ts j}[-1]\tss F_{2\ts i}[-1]-F_{1\ts j}[-1]\tss F_{1\ts i}[-1]\\
{}=F_{ij}[-2]-F_{ij}[-2]=0.
\non
\end{multline}
This implies that the terms in the expansion of $F_{12}[0]\ts\Pf F[-1]$
corresponding to pairs of the form $(\si,\si')$ cancel pairwise. Thus,
$F_{12}[0]\ts\Pf F[-1]=0$.

Now we verify that
\beql{fone}
F_{12}[1]\ts\Pf F[-1]=0.
\eeq
Consider first the summands in \eqref{pf} with
$\si(1)=1$ and $\si(2)=2$. In $V(\wh\oa_{2n})$ we have
\ben
F_{12}[1]\ts F_{1\ts 2}[-1]\tss F_{\si(3)\ts \si(4)}[-1]\dots F_{\si(2n-1)\ts\si(2n)}[-1]
=-K\ts F_{\si(3)\ts \si(4)}[-1]\dots F_{\si(2n-1)\ts\si(2n)}[-1].
\een
Furthermore, let $\tau\in\Bc_{2n}$
with $\tau(2)>2$. Then $\tau(3)=2$ and $\tau(4)>2$. We have
\ben
\bal
F_{12}[1]\ts F_{1\ts \tau(2)}[-1]\tss F_{2\ts \tau(4)}[-1]&\dots F_{\tau(2n-1)\ts\tau(2n)}[-1]\\
{}&=-F_{2\ts \tau(2)}[0]\tss F_{2\ts \tau(4)}[-1]\dots F_{\tau(2n-1)\ts\tau(2n)}[-1]\\
{}&= F_{\tau(2)\ts \tau(4)}[-1]\dots F_{\tau(2n-1)\ts\tau(2n)}[-1].
\eal
\een
Note that for any given $\si$, the number of elements $\tau$ such that the product
\ben
F_{\tau(2)\ts \tau(4)}[-1]\dots F_{\tau(2n-1)\ts\tau(2n)}[-1]
\een
coincides, up to a sign, with the product
\ben
F_{\si(3)\ts \si(4)}[-1]\dots F_{\si(2n-1)\ts\si(2n)}[-1]
\een
equals $2n-2$. Indeed, for each $k=2,\dots,n$ the unordered pair $\{\tau(2),\tau(4)\}$
can coincide with the pair $\{\si(2k-1),\si(2k)\}$. Hence, taking the signs
of permutations into account, we find that
\ben
F_{12}[1]\ts \Pf F[-1]
=(-K-2n+2)\ts \sum_{\si}\ts \sgn\si\cdot
F_{\si(3)\ts \si(4)}[-1]\dots F_{\si(2n-1)\ts\si(2n)}[-1],
\een
summed over $\si\in\Bc_{2n}$ with $\si(1)=1$ and $\si(2)=2$.
Since $-K-2n+2=0$ at the critical level, we get \eqref{fone}.
\epf

Introduce formal Laurent series
\ben
F_{ij}(z)=\sum_{r\in\ZZ}F_{ij}[r]\ts z^{-r-1}\Fand F_{ij}(z)_+=\sum_{r<0}F_{ij}[r]\ts z^{-r-1}
\een
and expand the Pfaffians of the matrices $F(z)=[F_{ij}(z)]$ and $F(z)_+=[F_{ij}(z)_+]$ by
\ben
\Pf F(z)=\sum_{p\in\ZZ} S_p\ts z^{-p-1}\Fand
\Pf F(z)_+=\sum_{p<0} S^+_p\ts z^{-p-1}.
\een
Invoking the vertex algebra structure on the vacuum module $V(\oa_{2n})$ (see \cite{f:lc}),
we derive from Theorem~\ref{thm:pfaff} that
the coefficients $S_p$ are {\it Sugawara operators} for $\wh\oa_{2n}$;
they commute with the elements of $\wh\oa_{2n}$ (note that normal ordering
is irrelevant here, as the coefficients of the series pairwise commute). Moreover,
the coefficients $S^+_p$ are elements of the Feigin--Frenkel center $\z(\wh\oa_{2n})$.


\end{document}